\documentclass[a4paper,12pt]{article}

\makeatletter
\@addtoreset{footnote}{page}
\makeatother

\usepackage{style-enumitem}

\usepackage{amsthm}
\usepackage{amsmath,amssymb,latexsym,amsfonts,mathrsfs}

\newcommand{\n}{\nonumber}

\renewcommand{\hat}{\widehat}

\renewcommand{\bar}{\overline}

\newcommand{\absol}[1]{\left| #1 \right|} 
\newcommand{\norm}[1]{\left\| #1 \right\|} 
\newcommand{\rbra}[1]{\!\left( #1 \right)} 
\newcommand{\cbra}[1]{\!\left\{ #1 \right\}} 
\newcommand{\sbra}[1]{\!\left[ #1 \right]} 

\newcommand{\bD}{\ensuremath{\mathbb{D}}}
\newcommand{\bE}{\ensuremath{\mathbb{E}}}

\newcommand{\bN}{\ensuremath{\mathbb{N}}}

\newcommand{\bP}{\ensuremath{\mathbb{P}}}

\newcommand{\bR}{\ensuremath{\mathbb{R}}}

\newcommand{\bT}{\ensuremath{\mathbb{T}}}

\theoremstyle{plain}
\newtheorem{Thm}{Theorem}[section]

\newtheorem{Lem}[Thm]{Lemma}
\newtheorem{Prop}[Thm]{Proposition}
\newtheorem{Cor}[Thm]{Corollary}

\theoremstyle{definition}

\newtheorem{Def}[Thm]{Definition}
\newtheorem{Rem}[Thm]{Remark}

\newcommand{\Proof}[2][Proof]{\begin{proof}[{#1}] #2 \end{proof}}

\numberwithin{equation}{section}

\makeatletter
\renewcommand\section{\@startsection {section}{1}{\z@}%
                                   {-3.5ex \@plus -1ex \@minus -.2ex}%
                                   {2.3ex \@plus.2ex}%
                                   {\normalfont\large\bf}}
\makeatother

\makeatletter
\renewcommand\subsection{\@startsection {subsection}{1}{\z@}%
                                   {-3.5ex \@plus -1ex \@minus -.2ex}%
                                   {2.3ex \@plus.2ex}%
                                   {\normalfont\normalsize\bf}}
\makeatother

\usepackage{color}
\renewcommand{\phi}{\hat{\psi}}

\allowdisplaybreaks[1]

\begin{document}

\begin{center}
{\Large \bf 
Approximation and duality problems \\of refracted processes
}
\end{center}
\begin{center}
Kei Noba
\end{center}

\begin{abstract}
For given two standard processes with no positive jumps, 
we construct, using the excursion theory, a Markov process whose positive and 
negative motions have the same law as the two processes. 
The resulting process is a generalization of Kyprianou--Loeffen's refracted L\'evy processes. We discuss approximation problem for our refracted processes coming from 
L\'evy processes by removing small jumps and taking the limit as the removal level tends to zero. We also discuss conditions for refracted processes to have dual processes.
\end{abstract}


\section{Introduction}
Let $X$ and $Y$ be $\bR$-valued standard processes with no positive jumps. 
We want to study 
a $\bR$-valued standard process $U$ whose positive and negative motions have the same law as X and Y, respectively. 
If such a process $U$ exists,  
we call $U$ a \textit{refracted process} (of $X$ and $Y$). 
We want to give its precise definition in a good generality 
and study approximation and duality problems for the process. 
\par
As an earlier result, Kyprianou--Loeffen \cite{KypLoe} 
studied a unique strong solution of the stochastic differential equation 
\begin{align}
U_t-U_0&=X_t-X_0 + \delta_U \int_0^t1_{\{U_s <0\}}ds,~~~t\geq 0, 
\end{align}
for a spectrally negative L\'evy process $X$ and 
a positive constant $\delta_U$. 
They called $U$ the \textit{refracted L\'evy process}. 
This process can be regarded in our terminology as a refracted process of $X$ and $Y$ 
with $Y_t\overset{d}{=} X_t + \delta_U t$. 
Noba--Yano \cite{NobYan} generalized Kyprianou--Loeffen's refracted L\'evy processes 
when $X$ and $Y$ have different Laplace exponents and $X$ has no 
Gaussian parts. 
When $X$ has unbounded variation paths and has no Gaussian parts, 
they used the excursion theory to construct refracted L\'evy processes 
from the law of stopped process and 
the excursion measure $n^U_0$ away from $0$ satisfying the following: 
for all non-negative measurable functional $F$, 
\begin{align}
n^U_0 \sbra{F(U) }
=
&n^X_0\sbra{\bE^{Y^0}_{X_{T^-_0}} \sbra{
F(w \circ Y^0 )}\big{|}_{w=k_{T^-_0}X}
;0< T^-_0\leq  T_0}     ,
\end{align}
where we denote by $n^X_0$ an excursion measure of $X$ away from $0$, 
by $Y^0$ the stopped process of $Y$ at $0$, 
by $T^-_0$ the first hitting time to $(-\infty , 0]$, by $T_0$ the first 
hitting time to $0$, by $\circ$ 
the concatenation of two c\`adl\`ag functions and 
by $k_{T^-_0}$ the killing operator at $T^-_0$. 
We always understand 
$\bE^{Y^0}_{X_{T^-_0}} \sbra{
F(w \circ Y^0 )}\big{|}_{w=k_{T^-_0}X}=F(X) $ on $\{  T^-_0 = T_0\}$. 
\par
In this paper, we use the excursion theory to construct a refracted process $U$ from given two standard processes $X$ and $Y$ with no positive jumps. 
For a non-negative constant $c_0$ and a negative function $\psi$, 
we define the excursion measure of $U$ away from $0$ 
in the following: 
for non-negative measurable functional $F$, 
\begin{align}
n^U_0 \sbra{F(U) }
&=c_0 n^Y_0\sbra{ F( Y) ;T^-_0=0} \notag
\\
&~+n^X_0\sbra{\bE^{Y^0}_{\psi(X_{T^-_0-}, X_{T^-_0})} \sbra{
F(w \circ Y^0 )}\big{|}_{w=k_{T^-_0}X}
;0< T^-_0\leq T_0}    . \label{304a}
\end{align}
{The resulting process $U$ constructed from $n^U_0$ via the excursion theory 
satisfies the following: 
\begin{align}
\text{For }& x<0, ~{\{U_t\}}_{t\leq T_0} ~(\text{under }\bP^U_x) 
\overset{d}{=}{\{Y_t\}}_{t\leq T_0} ~(\text{under }\bP^Y_x), \\
\text{For }& x>0, ~\rbra{{\{U_t\}}_{t < T^-_0}, U_{T^-_0}}~ (\text{under }\bP^U_x) 
\overset{d}{=}
\rbra{{\{X_t\}}_{t< T^-_0},\psi(X_{T^-_0-}, X_{T^-_0}) }~(\text{under }\bP^X_x).  
\end{align}}
We call $\psi$ the \textit{landing function} because it indicates the landing point at the 
first hitting time $T^-_0$ of $U$. 
\par
One of our main problems is approximation. 
We assume that our new refracted process $U$ comes from two L\'evy processes 
$X$ and $Y$. 
We will then prove that $U$ is the limit  in distribution on the c\`adl\`ag function 
space of the sequence ${\{U^{(n)}\}}_{n \in \bN}$ of refracted processes coming 
from the drifted compound Poisson processes constructed from $X$ and $Y$ by removing 
small jumps and by adding drifts. 
Noba--Yano \cite{NobYan} studied this problem in the special case of no Gaussian 
part of $X$, where our landing functions did not appear. 
In our setting, our landing functions play an important role: 
even when $U$ has a trivial landing function $\psi(x, y)=y$, the approximating 
process $U^{(n)}$ must involve a suitable landing function. 
\par
The other is duality. 
Let $X$ and $Y$ be standard processes with no positive jumps and 
let $\hat{X}$ and $\hat{Y}$ be dual processes of $X$ and $Y$, respectively. 
Let $U$ be the refracted process of $X$ and $Y$ and 
let $\hat{U}$ be the refracted process of $\hat{X}$ and $\hat{Y}$. 
We will then obtain the necessary and sufficient condition 
that the refracted processes $U$ and $\hat{U}$ are in duality 
in terms of a certain identity involving excursion measures 
and landing functions. 
To prove duality of $U$ and $\hat{U}$, we require that their excursion measures are 
transformed into each other by time reversal. 
For this purpose, we utilize landing functions in order to adapt the jumps at the 
switching time between $X$ and $Y$. 
\par
We give an example of a refracted process possessing a dual. 
We construct it from two spectrally negative stable processes, 
where we will make a computation to find a suitable landing function. 
\par
The organization of the present paper is as follows. 
In Section \ref{Pre} we propose some notation and recall preliminary facts about 
scale functions of standard processes with no positive jumps. 
In Section \ref{SecGen} we give the precise definition 
of our new refracted processes. 
In Section \ref{App} we study the approximation problem. 
In Section \ref{duality} we give the definition of duality. 
In Section \ref{SecDua} we study the duality problem. 
In Section \ref{Exa} we give an example of refracted processes 
in duality using stable processes.

\section{Preliminary}\label{Pre}
Let $\bR \cup \{\partial\}$ denote the one-point compactification of $\bR$. 
Let $\bD$ denote the set of functions 
$\omega: [0, \infty) \rightarrow \bR \cup \{ \partial \}$ which are c\`adl\`ag and satisfy 
\begin{align}
\omega(t)=\partial~~~~~~t \geq \zeta(\omega) 
\end{align}
where $\zeta (\omega)= \inf \{t>0 : \omega (t) = \partial \}$. 
Let $\cal{B}(\bD)$ denote the class of Borel sets of $\bD$ 
equipped with the Skorokhod topology. 
For $\omega \in \bD$, denote 
\begin{align}
&T^-_x (\omega):= \inf \cbra{ t > 0 : \omega(t) \leq x}, \\
&T^+_x (\omega):= \inf \cbra{ t > 0 : \omega (t ) \geq  x}, \\
&T_x (\omega):= \inf \cbra{t>0 : \omega (t) =x}. 
\end{align}
For $\omega, \omega_1, \omega_2 \in \bD$ and $s, t \in [0, \infty)$, 
we adopt the following notation:  
\begin{align}
&\rho_x \omega (t) 
= 
\begin{cases}
\omega(T_x - t -),~~~~~&t < T_x< \infty,
\\
x, ~~~~~~~~~~~~~~~~~~&t \geq T_x , 
\\
\partial, ~~~~~& t\geq 0, ~ T_x (\omega) = \infty,
\end{cases}
\\
&k_s \omega (t) 
= 
\begin{cases}
\omega( t ),~~~~~~~~~~~~~&t < s,
\\
\partial, ~~~~~~~~~~~~~~~~~&t \geq s , 
\end{cases} 
\\
&\omega_1 \circ \omega_2 (t)
=
\begin{cases}
\omega_1(t), ~~~~~~~~~~~~~~~~~~~~&t < \zeta ( \omega_1 ), \\
\omega_2(t- \zeta (\omega_1)),~~~~~~~~~~&t \geq \zeta (\omega_1),  
\end{cases}
\\
&\theta_s \omega (t) 
= \omega (t + s ).
\end{align}
We have introduced in \cite{Nob} 
the generalized scale functions of standard processes with no positive jumps. 
Let $\bT$ be an interval of $\bR$ and set $a_0 = \sup \bT$ and $b_0 =\inf \bT$. 
We assume that the process $(Z, \bP^Z_x)$ considered in this paper is a 
$\bT$-valued standard process with no positive jumps 
satisfying the following conditions: 
\begin{itemize}
\item[(A1)]{$(x, y) \mapsto \bE^Z_x\sbra{ e^{-T_y}}>0$ is a ${\cal{B}}(\bT) \times \cal{B}(\bT)$-measurable function.} 
\label{iv}
\item[(A2)]{$Z$ has a reference measure $m_Z$ on $\bT$, i.e. 
for $q \geq 0$ and $x\in \bT$, 
{the measure $R^{(q)}_Z1_{(\cdot)}(x)$ is absolutely continuous 
with respect to $m_Z(\cdot)$} 
where 
\begin{align}
R_Z^{(q)}f(x):= \bE^Z_x \sbra{ \int_0^\infty e^{-qt} f(Z_t) dt}
\end{align}
for non-negative measurable function $f$. 
Here and hereafter we use the notation $\int_b^a = \int_{(b, a]\cap\bR}$. 
In particular, $\int_{b-}^a=\int_{[b, a]\cap\bR}$. 
}
\end{itemize}
By \cite[Theorem $18.4$]{GemHor}, 
there exist a family of processes ${\{L^{Z, x}\}}_{x\in\bT}$ with 
$L^{Z, x}={\cbra{L^{Z,  x}_t}}_{t\geq 0}$ for 
$x \in \bT$ which we call \textit{local times} 
such that the following conditions hold: 
for all $q>0$, $x \in \bT$ and non-negative measurable function $f$ 
\begin{align}
&~~~\int_0^t f(Z_s)ds = \int_{\bT} f(y) L^{Z, y}_{t} m_Z (dy) ,   ~~~~\text{a.s.}  \label{201f}\\
&R_Z^{(q)}f(x) 
= \int_{\bT} f(y) \bE^Z_x\sbra{\int_{0}^\infty e^{-qt}dL^{Z, y}_{t} } m_Z (dy). \label{202f}
\end{align}
We have the following two cases: 
\begin{itemize}
\item{{\bf Case 1.} 
If $x\in \bT$ is regular for itself, 
this $L^{Z,x}$ is the continuous local time at $x$ \cite[pp.216]{BluGet}. 
Note that $L^{Z,x}$ has no ambiguity of multiple constant because of 
\eqref{201f} or \eqref{202f}. }
\item{{\bf Case 2.} 
If $x\in \bT$ is irregular for itself, 
we have 
\begin{align}
L^{Z, x}_t =l^Z_x \# \{0\leq s< t : Z_s=x \},~~~\text{a.s.}
\end{align}
for some constant $l^Z_x \in (0, \infty)$. 
}
\end{itemize}
\par
In Case 1, 
let $\eta^{Z, x}$ denote the inverse local time of $L^{Z, x}$. 
Let $n^Z_x$ be an excursion measure away from $x$ which is 
associated with $L^{Z, x}$ (See \cite{Ito}). 
Then, for all $q > 0$, we have 
\begin{align}
-\log \bE^Z_0\sbra{e^{-q\eta^{Z, x}(1)}}
=\delta^Z_x q + n^Z_x\sbra{1-e^{-qT_x} } 
\end{align}
for a non-negative constant $\delta^Z_x$ called the \textit{stagnancy rate}. 
We thus have 
\begin{align}
\bE^Z_x \sbra{\int_{0}^\infty e^{-qt} dL^{Z, x}_t} 
= \bE^Z_x\sbra{\int_0^\infty e^{-q\eta^{Z, x} (s)}ds} 
=\frac{1}{\delta^Z_x q +n^Z_x \sbra{1 - e^{-qT_x}}}. 
\label{203z} 
\end{align}
In Case 2, we define $n^Z_x = \frac{1}{l^Z_x}\bP^{Z^x}_x$ 
where $\bP^{Z^x}_x$ denotes the law of $Z$ started from $x$ and stopped at $x$. 
Then 
we have 
\begin{align}
&\bE^Z_x \sbra{ \int_{0-}^{\infty} e^{-qt} d L^{Z, x}_t }
=l^Z_x \sum_{ i = 0}^{\infty} \rbra{\bE^Z_x \sbra{e^{-qT_x}}}^i
=\frac{l^Z_x}{\bE^Z_x \sbra{ 1-e^{-qT_x}}}= \frac{1}{n^Z_x \sbra{1-e^{-qT_0} }}.  
\end{align}
In \cite[Definition 3.1]{Nob}, the author has introduced the $q$-scale function of $Z$ as, 
for $q\geq 0$ and $x, y \in \bT$,  
\begin{align}
W_Z^{(q)}(x, y) =
\begin{cases}
\frac{1}{n^Z_y \sbra{e^{-q T^+_x }; T^+_x < \infty }},  ~~~~~~~~~~   &x>y,
\\
0,~~~~~~~~~~~~~~~~~~~~~~~~~~~&x \leq y.
\end{cases}\label{216}
\end{align}
Let us fix $b, a \in \bT$ with $b<a$. 
We need the following results. 
\begin{Thm}[{\cite[Theorem 3.4]{Nob}}]\label{Lem102}
For $q\geq 0$ and $ x \in (b, a)$, we have 
\begin{align}
\bE^Z_x \sbra{ e^{-q T^+_a}; T^+_a < T^-_b}
= \frac{W_Z^{(q)}(x, b)}{W_Z^{(q)}(a , b) }.\label{103}
\end{align}
\end{Thm}
For $q \geq 0$, $x \in (b, a)$ and non-negative measurable function $f$, we define
\begin{align}
\overline{\underline{R}}_Z^{(q; b, a)}f(x)
:= \bE^Z_x\sbra{ \int_0^{T^-_b \land T^+_a} e^{-qt}f(Z_t)dt}.
\end{align}
Then, for $q \geq 0$, we have 
\begin{align}
\overline{\underline{R}}_Z^{(q; b, a)}f(x)
&= \int_\bT f(y) \bE^Z_{x} 
\sbra{\int_{0}^{T^-_b \land T^+_a} e^{-qt}dL^{Z, y}_{t} } m_Z (dy) .
\end{align}
\begin{Thm}[{\cite[Theorem 3.6]{Nob}}]\label{Lem105}
For $q \geq 0 $ and $x, y \in (b, a)$, we have
\begin{align}
\bE^Z_{x} \sbra{\int_{0}^{T^-_b \land T^+_a} e^{-qt}dL^{Z, y}_t }
= \frac{W_Z^{(q)}(x, b)}{W_Z^{(q)}(a, b)}W_Z^{(q)}(a , y) - W_Z^{(q)}(x, y). \label{114}
\end{align}
\end{Thm}
\begin{Lem}[{\cite[Lemma 3.5]{Nob}} and {\cite[Lemma 6.1]{NobYan}}]
	For $q \geq 0$ and $x\in (b, a)$, we have
		\begin{align}
		\bE^Z_{x}\sbra{ e^{-qT_{a}^{+}} ;T_{a}^{+} <T_{b}^{-}}
		&=n^Z_x\sbra{ e^{-qT_{a}^{+}}; T^+_a < \infty } 
\bE^Z_x \sbra{\int_{0-}^{T^+_a \land T^-_b}e^{-qt}dL^{Z, x}_t} \\
&=\frac{n^Z_x\sbra{ e^{-qT_{a}^{+}}; T^+_a < \infty } }
{\delta^Z_xq+n^Z_x\sbra{1-e^{-qT_x}1_{\{T^+_a=\infty , T^-_b=\infty\}} }}.
\label{122aa}
		\end{align}
\end{Lem}

\section{Refracted processes} \label{SecGen}
In this section, we construct a refracted process from 
two \bR-valued standard processes with 
no positive jumps $X$ and $Y$ using the excursion theory. 
\par
Let $a_0$, $a_1$, $b_0$ and $b_1$ be real numbers with 
$-\infty \leq b_0 \leq b_1 < 0 < a_1 \leq a_0 \leq \infty$. 
Let $\bT_X$ be an interval with $\sup\bT_X = a_0$ and $\inf\bT_X = b_1$. 
Let $\bT_Y$ be an interval with $\sup\bT_Y = a_1$ and $\inf\bT_X = b_0$. 
We let $\bT:=\bT_X \cup \bT_Y$.  
Let $X$ and $Y$ be $\bT_X$ and $\bT_Y$-valued standard processes 
with no positive jumps, respectively. 
We assume $X$ (resp. $Y$) satisfying the following conditions: 
\begin{itemize}
\item[(B1)]{$(x, y) \rightarrow \bE^X_x \sbra{e^{-T_y}}>0$ 
(resp. $(x, y) \rightarrow \bE^Y_x \sbra{e^{-T_y}}>0$) is a 
${\cal{B}}(\bT_X) \times {\cal{B}}(\bT_X)$
(resp. ${\cal{B}}(\bT_Y) \times {\cal{B}}(\bT_Y)$)-measurable function. }
\item[(B2)]
We assume that $\lim_{y\uparrow x}\bE^X_y\sbra{e^{-T_x}}=1$ for all 
$x\in \bT_X\cap(0, \infty)$
(resp. $\lim_{y\uparrow x}\bE^Y_y\sbra{e^{-T_x}}=1$ for all 
$x\in \bT_Y\cap (-\infty , 0]$). 
\item[(B3)]
If $a_0 \notin \bT_X$, we assume that $\lim_{x \uparrow a_0}\bE^X_x\sbra{ e^{- T^-_y}}=0$ for all $y\in\bT_X$
(resp. 
If $b_0 \notin \bT_Y$, we assume that $\lim_{x \downarrow b_0}\bE^Y_x\sbra{ e^{- T^+_y}}=0$ for all $y\in\bT_Y$). 
\item[(B4)]{$X$ (resp. $Y$) has a reference measure $m_X$ on $\bT_X$ 
(resp. $m_Y$ on $\bT_Y$).  }
\end{itemize}
We define local times ${\{L^{X, x}\}}_{x \in \bT_X}$ and ${\{L^{Y, x}\}}_{x \in \bT_Y}$, 
excursion measures ${\{n^X_x\}}_{x\in\bT_X}$ and ${\{n^Y_x\}}_{x\in\bT_Y}$, 
and scale functions ${\{W_X^{(q)}\}}_{q \geq 0}$ and ${\{W_Y^{(q)}\}}_{q \geq 0}$ 
of $X$ and $Y$ in the same way as $Z$'s in Section \ref{Pre}, respectively. 
\par
Let
$\psi: (0, \infty)\times (-\infty ,  0) \rightarrow (- \infty , 0)$ 
be a measurable function satisfying 
\begin{align}
n^X_0 \sbra{1- e^{-T^-_0} \bE^Y_{
J_X } 
\sbra{e^{-T_0}} ;0<T^-_0 < T_0} < \infty, \label{301i}
\end{align}
where $J_X=\psi(X_{ T^-_0-} , X_{ T^-_0})$. 
Let $c_0 \geq 0$ be a constant. 
We define the law of stopped process $\bP^{U^0}_x$ 
for $x \in\bT \backslash \{0\}$ and 
the excursion measure $n^U_0$ away from $0$ by 
\begin{align}
\bE^{U^0}_x\sbra{F(U^0) }
&=
\begin{cases}
\bE^{Y^0}_x \sbra{f(Y^0)}, ~~~~~~~~~~~~~~~~~~~~~~~~~~~~~~~~~~~~~~~~~~~~~~~~&x \in \bT \cap(-\infty , 0),
\\
\bE^X_x\sbra{\bE^{Y^0}_{J_X} \sbra{
F(w \circ Y^0 )}\big{|}_{w=k_{T^-_0} X}; T^-_0 \leq T_0}
, ~&x\in \bT \cap (0, \infty),
\end{cases}
\label{302i}
\\
n^U_0 \sbra{F(U) }
&=c_0 n^Y_0\sbra{ F( Y) ;T^-_0=0} \notag
\\
&~+n^X_0\sbra{\bE^{Y^0}_{J_X} \sbra{
F(w \circ Y^0 )}\big{|}_{w=k_{T^-_0}X}
;0< T^-_0\leq T_0}     \tag{\ref{304a}}
\end{align}
for all non-negative measurable functional $F$ 
(if $\bP^X_0\sbra{T_0>0}=1$ or $\bP^Y_0\sbra{T_0>0}=1$, we assume that $c_0 = 0$). 
We write  $X^0$ and $Y^0$ for the stopped processes of $X$ and $Y$ 
upon hitting zero, respectively. 
By means of the excursion theory, we can construct 
from $n^U_0$ and ${\{\bP^{U^0}_x\}}_{x\in \bT\backslash \{0\}}$ 
a $\bT$-valued right continuous strong Markov process 
without stagnancy at $0$ (See, e.g., \cite{Sal}). 
\begin{Rem}
The condition $c_0 = 0$ is necessary when $\bP^X_0\sbra{T_0>0}=1$. 
Indeed, when $\bP^X_0\sbra{T_0>0}=1$ and $c_0 > 0$, 
the measure $n^U_0$ does not satisfy 
the condition \cite[pp.323, (vi')]{Sal} 
and  then $n^U_0$ is not an excursion measure.  
\end{Rem}
\begin{Lem}
The refracted process $U$ is a Feller process. So $U$ is a standard process. 
\end{Lem}
\Proof{
Let $C_0(= C^{\bT}_0)$ 
denote the set of continuous functions $f$ from $\bT$ to $\bR$ 
such that $f(x)\rightarrow 0$ as $x\downarrow b_0 $ when $b_0 \notin \bT$ 
and as $x\uparrow a_0 $ when $a_0 \notin \bT$. 
For $f \in C_0$, we write $\norm{f}=\sup_{x \in \bR} \absol{f(x)}$.
It is sufficient to verify the following conditions:
\begin{enumerate}
\item{For all $q > 0$, $R_U^{(q)}$ is a map from $C_0$ to $C_0$.}\label{3i}
\item{For all $f \in C_0$, $\lim_{q\uparrow \infty} \norm{ qR_U^{(q)} f - f}=0$. }
\label{3ii}
\end{enumerate}
$1)$ {\bf{The proof of {\eqref{3i}}}}\par
First, we prove that $R_U^{(q)}f$ is continuous. 
We let $x \in \bT$. 
By the construction of $U$ and (B2), it is easy to check that 
$\lim_{y \uparrow x}\bE^U_y \sbra{e^{-qT_x} }
=\lim_{y \downarrow x}\bE^U_x \sbra{e^{-qT_y} }=1$. 
We fix $x \in \bT$. 
For $y<x$, we have
\begin{align}
&\bar{\lim}_{y \uparrow x}\absol{ R_U^{(q)}f(x ) - R_U^{(q)}f(y )}\\
\leq& 
\bar{\lim}_{y \uparrow x}\absol{ R_U^{(q)}f(x ) - \bE^U_y\sbra{e^{-qT_x} }R_U^{(q)}f(x )}
+\bar{\lim}_{y \uparrow x}\absol{ \bE^U_y\sbra{\int_0^{T_x} e^{-qt} f(U_t) dt} }=0. 
\end{align}
For $y>x$, we have
\begin{align}
&\bar{\lim}_{y \downarrow x}\absol{ R_U^{(q)}f(x ) - R_U^{(q)}f(y )} \\
\leq &\bar{\lim}_{y \downarrow x} \absol{ 
\bE^U_x \sbra{e^{-q T_y} }R_U^{(q)}f(y )-R_U^{(q)}f(y )}
+\bar{\lim}_{y \downarrow x} \absol{\bE^U_x\sbra{\int_0^{T_y} e^{-qt} f(U_t) dt}}=0.
\end{align}
Second, we prove that  
$\lim_{x\uparrow a_0} R_U^{(q)}f(x)=0$ when $a_0 \notin \bT$ 
and $\lim_{x\downarrow b_0} R_U^{(q)}f(x)=0$ when $b_0\notin \bT$. 
We assume that $a_0 \notin \bT$. 
By {the assumption (B3)}, for all $x \in (0, a_0)$, 
$\lim_{y \uparrow a_0} \bE^U_y\sbra{ e^{-T^-_x}} 
= \lim_{y \uparrow a_0}\bE^X_y\sbra{ e^{-T^-_{x}}}=0$. 
Since $f \in C_0$, for all $\epsilon >0$, there exists $\delta \in (0, a_0)$ 
such that $\sup_{x \in (\delta , a_0)} \absol{f(x)} < \epsilon$.  
So we have
\begin{align}
\lim_{x \uparrow a_0}\absol{R_U^{(q)}f(x)}
&\leq \lim_{x \uparrow a_0} \rbra{
\bE^X_x \sbra{\int_0^{T^-_\delta} e^{-qt}\absol{f(X_t)} dt }
+
\bE^U_x \sbra{\int_{T^-_\delta}^\infty e^{-qt} \norm{f} dt }}
\\
&\leq  \frac{\epsilon}{q}
+\lim_{x \uparrow a_0} \bE^X_x\sbra{e^{-qT^-_\delta}} 
\frac{\norm{f}}{q} =\frac{\epsilon}{q}.
\end{align}
Therefore we have $\lim_{x \uparrow a_0}\absol{R_U^{(q)}f(x)}=0$. 
In the same way, we have 
$\lim_{x \downarrow b_0}\absol{R_U^{(q)}f(x)}=0$ when $b_0\notin \bT$. 
\par
$2)$ {\bf{The proof of {\eqref{3ii}}}}\par
By classical arguments, it is sufficient to prove 
$\lim_{q\uparrow \infty} \absol{qR^{(q)}_Uf(x)-f(x)} =0 $ 
for $x \in \bT$. Fix $x\in \bT$. 
For all $\epsilon >0$, there exists $\delta>0$ such that 
\begin{align}
\absol{x - y} < \delta \Rightarrow \absol{f(x) - f(y)} < \epsilon,~~~~~~~
x, y\in \bT.
\end{align}
We define 
\begin{align}
T^{\uparrow}_\delta = \inf \cbra{t> 0 : \absol{U_t -x} \geq \delta   }. 
\end{align}
Then we have 
\begin{align}
\absol{qR^{(q)}_U f(x) - f(x)}
&\leq q\bE^U_x\sbra{\int_0^{T^{\uparrow}_\delta}e^{-qt} \absol{f(U_t)-f(x)} dt} 
+q\bE^U_x\sbra{\int_{T^{\uparrow}_\delta}^\infty e^{-qt} \absol{f(U_t)-f(x)} dt} \\
&\leq \epsilon \bE^U_x \sbra{1 - e^{-qT^{\uparrow}_\delta} }
+2 \norm{f}\bE^U_x \sbra{ e^{-qT^{\uparrow}_\delta} }. \label{315z}
\end{align}
By the dominated convergence theorem, we have 
\begin{align}
\limsup_{q \uparrow \infty} \absol{qR^{(q)}_U f(x) - f(x)}\leq \epsilon. 
\end{align}
and so we have $\lim_{q \uparrow \infty}\absol{qR^{(q)}_U f(x) - f(x)}=0$. 
The proof is completed. 
}

\section{Approximation problem for refracted processes 
coming from L\'evy processes}\label{App}
Let $U$ be the refracted process 
constructed by $X$, $Y$, $\psi$ and $c_0$ as 
Section \ref{SecGen}. 
In this section, 
we assume that $X$ and $Y$ are spectrally negative L\'evy processes and 
we shall construct a sequence ${\{U^{(n)}\}}_{n\in\bN}$ of refracted processes 
coming from compound Poisson processes 
which converges to $U$ in distribution. 
In \cite[Section 8]{NobYan}, Noba--Yano studied this approximation problem 
only when $\sigma_X=0 $ and $\psi (x, y)=y$. 
So this section is a generalization of \cite[Section 8]{NobYan}. 
\begin{itemize}
\item[(C0)]
Let $X$, $Y$ be spectrally negative L\'evy processes which have
Laplace transforms 
\begin{align}
&\Psi_X(\lambda ) = 
\chi_X \lambda + \frac{\sigma^2_X}{2} \lambda^2 
-\int_{(-\infty , 0)} (1 - e^{\lambda y} + \lambda y1_{(-1 , 0)}(y))\Pi_X (dy), ~~~\lambda \geq 0, \\
&\Psi_Y(\lambda ) = 
\chi_Y \lambda + \frac{\sigma^2_Y}{2} \lambda^2 
-\int_{(-\infty , 0)} (1 - e^{\lambda y} + \lambda y1_{(-1 , 0)}(y))\Pi_Y (dy),~~~\lambda \geq 0
\end{align}
for some constants $\chi_X , \chi_Y \in \bR$, $\sigma_X , \sigma_Y \geq 0$
and some L\'evy measures $\Pi_X, \Pi_Y$, respectively. 
We let $\Phi_X(\theta)=\inf \{\lambda>0 : \Psi_X (\lambda)> \theta\} $ 
and $\Phi_Y (\theta)=\inf \{\lambda>0 : \Psi_Y (\lambda)> \theta\} $.
We assume that reference measures $m_X$, $m_Y$ are Lebesgue measures 
and let the excursion measures ${\{n^X_x\}}_{x \in \bR}$ and ${\{n^Y_x\}}_{x \in \bR}$ 
of $X$ and $Y$ be those in Section \ref{SecGen} 
satisfying the following, respectively: 
for $x\in\bR$ and $q > 0$, 
\begin{align}
n^X_x \sbra{1-e^{-qT_x} }
=\frac{1}{\Phi_X^{\prime}(q)}, ~~~
n^Y_x \sbra{1-e^{-qT_x} }
=\frac{1}{\Phi_Y^{\prime}(q)}. 
\end{align}
Let $\psi$ be a continuous landing function 
which has the following condition: 
\begin{align}
\text{There exist }k, l>0 \text{ such that }
\psi (x, y) \geq l (y-x),\text{ for } x-y<k . \label{404}
\end{align}
{(Note that \eqref{404} implies \eqref{301i}.)} 
Let $c_0 $ be a non-negative constant such that $c_0=0$ when $\sigma_X=0$ 
or $\sigma_Y=0$. 
\item[(C1)]
Let ${\{\epsilon^X_n \}}_{n\in\bN}$ and 
${\{\epsilon^Y_n\}}_{n \in \bN}$ be sequences of strictly positive numbers 
satisfying  
\begin{align}
&\lim_{n\uparrow \infty}\epsilon^X_n =\lim_{n\uparrow \infty}\epsilon^Y_n=0. 
\end{align}
When $c_0>0$ (and consequently $\sigma_X\sigma_Y>0$), we assume that 
\begin{align}
&\lim_{n\uparrow \infty} \frac{\epsilon^Y_n}{\epsilon^X_n}
=\frac{\sigma_Y^2}{\sigma_X^2}c_0.\label{403ab}
\end{align}
For $n\in \bN$, we define 
		\begin{align}
		\Psi_{X^{(n)}} (\lambda) &=  \chi_X \lambda  
		-\frac{\sigma_{X}^{2}}{{(\epsilon^{X}_n)}^2} \rbra{ 1 - e^{ \lambda(-\epsilon^X_n)} 
			+ \lambda \rbra{ -\epsilon^X_n} } \nonumber\\
		&~~~~~~~~- \int_{\rbra{-\infty , -\epsilon^X_n}}\rbra{ 1 - e^{\lambda y} + 
		\lambda y 1_{\rbra{-1, -\epsilon^X_n}} (y)} \Pi_X (dy)
		\\
		&=\delta_{X^{(n)}}\lambda-\int_{(-\infty, 0)}\rbra{1-e^{\lambda y}}\Pi_{X^{(n)}}(dy)
		\end{align}
	where 
	\begin{align}
	\delta_{X^{(n)}}&=\chi_X+\frac{\sigma_{X}^2 }{\epsilon^X_n} + \int_{(-1, -\epsilon^X_n)}(-y)\Pi_X(dy)
	\\
	\Pi_{X^{(n)}} &=1_{(-\infty, -\epsilon^X_n)}\Pi_X 
	+\frac{\sigma_X^2}{{(\epsilon^{X}_n)}^2} \delta_{(-\epsilon^X_n)}.
	\end{align}
	Let $X^{(n)}$ be a compound Poisson process with positive drift
	which has Laplace exponent $\Psi_{X^{(n)}}$. 
	We let $\Phi_{X^{(n)}}$ denote the right inverse of $\Psi_{X^{(n)}}$. 
	We note that $\Psi_{X^{(n)}}(\lambda) \rightarrow \Psi_X(\lambda)$ for all $\lambda \geq 0$, so 
	that we have $X^{(n)} \rightarrow X$ in law on $\bD$. 
{Furthermore, we have $\lim_{n\uparrow \infty}\Phi_{X^{(n)}}(\lambda) = \Phi_X(\lambda)$ for all $\lambda \geq 0$.} 
	More preciously, by \cite[pp.210]{Ber}, we see that there exists a coupling 
	of $X^{(n)}$'s such that $X^{(n)} \rightarrow X$ uniformly on compact intervals 
	almost surely. 
	We define $\Psi_{Y^{(n)}}$, $\delta_{Y^{(n)}}$, $\Pi_{Y^{(n)}}$, $\Phi_{Y^{(n)}}$ and $Y^{(n)}$ 
	in the same way as those for $X$. 
\end{itemize}
		\begin{Lem}\label{Lem402y}
We assume that $\sigma_Y>0$. 
	Then for all $q >0$ and all bounded continuous function $f$, we have 
		\begin{align}
		\lim_{n \uparrow \infty} \frac{\sigma_Y^2}{{(\epsilon^{Y}_n)}^2} 
\int_{0}^{\epsilon^Y_n} 
			R_{k_{T^+_0}Y^{(n)}}^{(q)} f\rbra{-v}dv
		=n_0^Y \sbra{ \int_0^{T_0} e^{-qt}f(Y_t )dt ; T^-_0 = 0}. 
		\label{2-503}
		\end{align}
	\end{Lem}
\Proof{
By the definition of ${\{Y^{(n)}\}}_{n\in\bN}$, 
we have that for all $q>0$, 
\begin{align}
\lim_{n \uparrow \infty} \frac{n^{Y^{(n)}}_0\sbra{\int_0^\infty e^{-qt} g(Y^{(n)}_t) dt}}
{n^{Y^{(n)}}_0\sbra{\int_0^\infty e^{-qt} dt}}
=\lim_{n \uparrow \infty} R^{(q)}_{Y^{(n)}}g(0)=R^{(q)}_Y g(0)
=\frac{n^{Y}_0\sbra{\int_0^\infty e^{-qt} g(Y_t) dt}}
{n^{Y}_0\sbra{\int_0^\infty e^{-qt} dt}} \label{411}
\end{align}
and $\lim_{n\uparrow \infty} R^{(q)}_{Y^{(n)0}}f(u) =R^{(q)}_{Y^{0}}f(u)$ 
for $u<0$ and for $g=f1_{(-\infty , 0)}$ or $f1_{(0, \infty)}$.  
By \cite[Lemma 3.5]{NobYan} (which can easily be extended to the case of positive Gaussian component) and 
by $\lim_{n\uparrow \infty}\Phi_{Y^{(n)}}(\lambda) = \Phi_Y (\lambda)$ on 
for all $\lambda \geq 0$, we have, for all $q> 0$, 
\begin{align}
\lim_{n \uparrow \infty}
n^{Y^{(n)}}_0\sbra{\int_0^\infty e^{-qt} f(Y^{(n)}_t)1_{(0, \infty)}(Y^{(n)}_t) dt}
=n^{Y}_0\sbra{\int_0^\infty e^{-qt} f(Y_t)1_{(0 , \infty)}(Y_t) dt}. 
\end{align}
and thus by \eqref{411}, 
we have $\lim_{n \uparrow \infty}n^{Y^{(n)}}_0\sbra{\int_0^\infty e^{-qt} dt}
=n^{Y}_0\sbra{\int_0^\infty e^{-qt} dt}$. 
Again by \eqref{411}, we obtain 
\begin{align}
\lim_{n \uparrow \infty}
n^{Y^{(n)}}_0\sbra{\int_0^\infty e^{-qt} f(Y^{(n)}_t)1_{(-\infty , 0)}(Y^{(n)}_t)dt}
=n^{Y}_0\sbra{\int_0^\infty e^{-qt} f(Y_t)1_{(-\infty , 0)}(Y_t)dt}.  \label{412y}
\end{align}
By \cite[Theorem 3.3]{NobYan} 
(which can easily be extended to the case of positive Gaussian component), we have
\begin{align}
&n^Y_0 \sbra{\int_0^{T_0} e^{-qt}f(Y_t)1_{(-\infty, 0)}(Y_t) dt}\\
=&n^Y_0 \sbra{\int_0^{T_0} e^{-qt}f(Y_t) dt ; T^-_0 = 0}
+n^Y_0 \sbra{e^{-qT^-_0}\bE^{Y}_{Y_{T^-_0}} \sbra{\int_0^{T_0}e^{-qt} f(Y_t)dt }
; 0<T^-_0 < T_0} \\
=&n^Y_0 \sbra{\int_0^{T_0} e^{-qt}f(Y_t) dt ; T^-_0 = 0}
+\int_0^\infty dv \int_{(-\infty , 0)} R_{k_{T^+_0}Y}^{(q)}f(u)
e^{-\Phi_Y (q) v} \Pi_Y (du - v). \label{415y}
\end{align}
and 
\begin{align}
&n^{Y^{(n)}}_0 \sbra{ \int_0^{T_0} e^{-qt}f(Y^{(n)}_t)1_{(-\infty,  0)}(Y^{(n)}_t)dt} \n\\ 
=&n^{Y^{(n)}}_0 \sbra{e^{-qT^-_0}\bE^{Y^{(n)}}_{Y^{(n)}_{T^-_0}} \sbra{\int_0^{T_0}e^{-qt} f(Y^{(n)}_t)dt }
; T^-_0 < T_0} \\
=&\int_0^\infty dv\int_{(-\infty , 0)} R_{k_{T^+_0}Y^{(n)}}^{(q)}f(u)
e^{-\Phi_{Y^{(n)}} (q) v} \Pi_{Y^{(n)}} (du - v)\\
=&\frac{\sigma_Y^2 }{{(\epsilon^{Y}_n)}^2}
\int_0^{\epsilon^{Y}_n}R^{(q)}_{k_{T^+_0}Y^{(n)}}f\rbra{v- \epsilon^{Y}_n}
e^{-\Phi_{Y^{(n)}}(q)v}dv \n \\
&+\int_0^\infty dv \int_{(-\infty , 0\land(-\epsilon^{Y}_n+v))}
 R_{k_{T^+_0}Y^{(n)}}^{(q)}f(u)
e^{-\Phi_{Y^{(n)}} (q) v} \Pi_{Y} (du - v). \label{418y}
\end{align}
By the same argument as that of the proof of \cite[Theorem 8.4]{NobYan}, we have 
\begin{align}
&\lim_{n\uparrow \infty}
\int_0^\infty dv\int_{(-\infty , 0\land(-\epsilon^{Y}_n+v))}
 R_{k_{T^+_0}Y^{(n)}}^{(q)}f(u)
e^{-\Phi_{Y^{(n)}} (q) v} \Pi_{Y} (du - v) \n\\
=&\int_0^\infty dv\int_{(-\infty , 0)} R_{k_{T^+_0}Y}^{(q)}f(u)
e^{-\Phi_Y (q) v} \Pi_Y (du - v) \label{419y}
\end{align}
By \eqref{412y}, \eqref{415y}, \eqref{418y} and \eqref{419y}, we obtain
\begin{align}
\lim_{n \uparrow \infty}\frac{\sigma_Y^2}{{(\epsilon^{Y}_n)}^2} 
 \int_{0}^{\epsilon^{Y}_n} 
			R_{k_{T^+_0}Y^{(n)}}^{(q)} f\rbra{v-\epsilon^{Y}_n}e^{-\Phi_{Y^{(n)}}(q)v}dv
		=&n_0^Y \sbra{ \int_0^{T_0} e^{-qt}f(Y_t )dt ; T^-_0 = 0}. \label{422}
\end{align}
By a simple argument, we can see that the left hand side of \eqref{422} coincides with 
that of \eqref{2-503}, which leads to the desired conclusion. 
}
\begin{itemize}
\item[(C2)]
Let ${\{\psi^{(n)}\}}_{n \in \bN}$ be a sequence of functions satisfying 
\begin{align}
\psi^{(n)}(x, y)= \psi(x, y)1_{\{ x-y >\epsilon^X_n \}} 
-\frac{\sigma_Y^2}{\sigma_X^2}c_0 x 1_{\{x - y=\epsilon^X_n\}}
\end{align} 
for all $x>0$, $y< 0$ and $n \in \bN$ where we understand $\frac{0}{0}=0$. 
\end{itemize}
\begin{Thm}\label{Thm404}
Let $X^{(n)}$ and $Y^{(n)}$ be those in (C1) and let $\psi^{(n)}$ be that in (C2). 
Let $U^{(n)}$ be the refracted process 
constructed by $X^{(n)}$, $Y^{(n)}$, $\psi^{(n)}$ and $c_0^{(n)}=0$. 
Then, for all $q > 0$, $x \in \bR$ and bounded continuous function $f$, we have 
\begin{align}
\lim_{n\uparrow \infty}R_{U^{(n)}}^{(q)}f(x) =  R_U^{(q)}f(x) .\label{424a}
\end{align}
\end{Thm}
\Proof{ 
$i$) We prove \eqref{424a} for $x=0$. 
For this purpose we shall prove that
\begin{align}
\lim_{n\uparrow \infty}
n^{U^{(n)}}_0 \sbra{ \int_0^{T_0} e^{-qt} f(U^{(n)}_t) dt}
=n^{U}_0 \sbra{ \int_0^{T_0} e^{-qt} f(U_t) dt}. \label{423ab}
\end{align}
for all $q>0$ and bounded continuous function $f$. 
By \cite[Lemma 3.5]{NobYan} and 
$\lim_{n \uparrow \infty }\Phi_{X^{(n)}} (\lambda)= \Phi_X (\lambda)$ 
for all $\lambda \geq 0$, 
we have 
\begin{align}
\lim_{n\uparrow \infty}
n^{U^{(n)}}_0 \sbra{ \int_0^{T^-_0} e^{-qt} f(U^{(n)}_t)  dt}
&=\lim_{n\uparrow \infty}
n^{X^{(n)}}_0 \sbra{ \int_0^{T^-_0} e^{-qt} f(X^{(n)}_t)  dt} \\
&=n^{X}_0 \sbra{ \int_0^{T^-_0} e^{-qt} f(X_t) dt}\\
&=n^{U}_0 \sbra{ \int_0^{T^-_0} e^{-qt} f(U_t) dt} .\label{433z}
\end{align}
By the definition of $n^{U^{(n)}}_0$ and \cite[Theorem 3.3]{NobYan}, 
we have
\begin{align}
&
n^{U^{(n)}}_0 
\sbra{ \int_{T^-_0}^{T_0} e^{-qt} f(U^{(n)}_t)  dt}\n\\
=&
n^{X^{(n)}}_0 \sbra{ e^{-qT^-_0}\bE^{Y^{(n)}}_{\psi^{(n)} (X_{T^-_0 -} , X_{T^-_0})}
\sbra{\int_0^{T_0} e^{-qt} f(U^{(n)}_t)  dt}; T^-_0 < T_0} \\
=&\int_0^\infty dv\int_{(-\infty , 0)} R_{k_{T^+_0}Y^{(n)}}^{(q)}f(\psi^{(n)} (v, u))
e^{-\Phi_{X^{(n)}} (q) v} \Pi_{X^{(n)}} (du - v)\\
=&\frac{\sigma_X^2}{{(\epsilon^{X}_n)}^2}
\int_0^{\epsilon^X_n}R^{(q)}_{k_{T^+_0}Y^{(n)}}
f\rbra{\psi^{(n)} (v, v- \epsilon^X_n)}
e^{-\Phi_{X^{(n)}}(q)v}dv \n\\
&~~~~~+\int_0^\infty dv \int_{(-\infty , 0\land(-\epsilon^X_n+v))} R_{k_{T^+_0}Y^{(n)}}^{(q)}f(\psi (v,  u))
e^{-\Phi_{X^{(n)}} (q) v} \Pi_{X} (du - v) \\
=&(\text{I}) + (\text{II}). \label{429ab}
\end{align}
\par
{Let us compute the limit of (II). 
We have 
\begin{align}
(\text{II})
=
\int_{(-\infty , 0)}\Pi_X (du) 1_{\{u < -\epsilon^Y_n\}}
\int_0^{-u} R_{k_{T^+_0}Y^{(n)}}^{(q)}f(\psi (v,  u+v))e^{-\Phi_{X^{(n)}} (q) v}dv \label{430}
\end{align}
To use the dominated convergence theorem, we dominate the integrand as 
\begin{align}
&\absol{1_{\{u < -\epsilon^Y_n\}} \int_0^{-u} R_{k_{T^+_0}Y^{(n)}}^{(q)}
f(\psi (v,  u+v))e^{-\Phi_{X^{(n)}} (q) v} dv}\\
\leq & 
\norm{f}\int_0^{-u} e^{- \Phi^{\inf}_X(q) v}\bE^{Y^{(n)}}_{\psi (v,  u+v)}
\sbra{\int_0^{T_0}e^{-qt}dt}dv \\
=&\frac{\norm{f}}{q} \int_0^{-u} e^{- \Phi^{\inf}_X(q) v}
(1-e^{\Phi^{\inf}_Y(q)\psi (v,  u+v)})dv, \label{433}
\end{align}
where $\Phi^{\inf}_X(q)=\inf_{n\in\bN}\Phi_{X^{(n)}}(q)$ and 
$\Phi^{\inf}_Y(q)=\inf_{n\in\bN}\Phi_{Y^{(n)}}(q)$. 
By \eqref{404}, we have 
\begin{align}
\eqref{433}&\leq
\frac{\norm{f}}{q}\int_0^{-u}e^{- \Phi^{\inf}_X(q) v}
(1-1_{\{u>-k\}}e^{\Phi^{\inf}_Y(q)l u})dv \\
&=\frac{\norm{f}}{q \Phi^{\inf}_X(q)}
(1-e^{\Phi^{\inf}_X(q) u})(1-1_{\{u>-k\}}e^{\Phi^{\inf}_Y(q)l u})\in L^1(\Pi_X). 
\label{435}
\end{align}
By \eqref{435} and the dominated convergence theorem, we have 
\begin{align}
\lim_{n\uparrow  \infty}\eqref{430}
=&\int_{(-\infty , 0)}\Pi_X (du) 
\int_0^{-u} R_{k_{T^+_0}Y}^{(q)}f(\psi (v,  u+v))e^{-\Phi_{X} (q) v}dv\\
=&\int_0^\infty dv\int_{(-\infty , 0)} R_{k_{T^+_0}Y}^{(q)}f(\psi(v, u))
e^{-\Phi_{X} (q) v} \Pi_{X} (du - v). \label{437}
\end{align}
By the definition of $n^U_0$ and \cite[Theorem 3.3]{NobYan}, 
we have
\begin{align}
\eqref{437}=
&n^X_0\sbra{e^{-qT^-_0} \bE^Y_{\psi(X_{T^-_0-}, X_{T^-_0})} 
\sbra{\int_0^{T_0} e^{-qt}f(Y_t) dt}; 0< T^-_0 < T_0}
\\
=&n^U_0\sbra{\int_{T^-_0}^{T_0} e^{-qt}f(U_t) dt ; 0< T^-_0 < T_0}. \label{431ab}
\end{align}
\par
Let us compute the limit of (I). }
Let $c_1 = \frac{\sigma_Y^2}{\sigma_X^2}c_0$.
By the definition of $\psi^{(n)}$, we have 
\begin{align}
(\text{I})
=
\frac{\sigma_X^2}{{(\epsilon^{X}_n)}^2}
\int_0^{\epsilon^X_n}R^{(q)}_{k_{T^+_0}Y^{(n)}}
f\rbra{-c_1v}e^{-\Phi_{X^{(n)}}(q)v}dv. \label{432ac}
\end{align}
When $c_1=0$, we have $\eqref{432ac}=0$. 
When $c_1>0$, we have 
\begin{align}
\lim_{n\uparrow\infty} \eqref{432ac}
=& \lim_{n\uparrow\infty} c_0c_1\frac{\sigma_Y^2}{{(\epsilon^{Y}_n)}^2}
\int_0^{\epsilon^X_n}R^{(q)}_{k_{T^+_0}Y^{(n)}}f\rbra{-c_1v }dv. \label{433ab}
\end{align}
By the change of variables, 
we have  
\begin{align}
c_0c_1\frac{\sigma_Y^2}{{(\epsilon^{Y}_n)}^2}
\int_0^{\epsilon^X_n}R^{(q)}_{k_{T^+_0}Y^{(n)}}f\rbra{-c_1v }dv=
c_0\frac{\sigma_Y^2}{{(\epsilon^{Y}_n)}^2}
\int_0^{c_1 \epsilon^X_n}R^{(q)}_{k_{T^+_0}Y^{(n)}}f\rbra{-v }dv. \label{434ac}
\end{align}
We prove 
\begin{align}
\lim_{n\uparrow \infty}c_0\frac{\sigma_Y^2}{{(\epsilon^{Y}_n)}^2}
\int_0^{c_1 \epsilon^X_n}R^{(q)}_{k_{T^+_0}Y^{(n)}}f\rbra{-v }dv
=\lim_{n\uparrow \infty}c_0\frac{\sigma_Y^2}{{(\epsilon^{Y}_n)}^2}
\int_0^{\epsilon^Y_n}R^{(q)}_{k_{T^+_0}Y^{(n)}}f\rbra{-v }dv. \label{434}
\end{align}
Let $M_Y (q)= \sup_{n\in\bN} \Phi_{Y^{(n)}}(q) \times (1\lor \sup_{n\in\bN} 
\frac{c_1 \epsilon^X_n}{\epsilon^Y_n})$. 
We have 
\begin{align}
&\absol{c_0\frac{\sigma_Y^2}{{(\epsilon^{Y}_n)}^2}
\int_0^{c_1 \epsilon^X_n}R^{(q)}_{k_{T^+_0}Y^{(n)}}f\rbra{-v }dv 
-c_0\frac{\sigma_Y^2}{{(\epsilon^{Y}_n)}^2}
\int_0^{\epsilon^Y_n}R^{(q)}_{k_{T^+_0}Y^{(n)}}f\rbra{-v }dv}\\
&\leq c_0\frac{\sigma_Y^2}{{(\epsilon^{Y}_n)}^2}\absol{c_1\epsilon^X_n - \epsilon^Y_n}
\norm{f}\sup_{0\leq v \leq (c_1 \epsilon^X_n) \lor \epsilon^Y_n}
\bE^{Y^{(n)}}_{-v} \sbra{\int_0^{T_0} e^{-qt}dt}
\\
&\leq\frac{c_0 \sigma_Y^2 \norm{f}}{ q}
\absol{c_1\frac{\epsilon^X_n}{\epsilon^Y_n} - 1} 
\frac{1-e^{-M_Y(q) \epsilon^Y_n}}{\epsilon^Y_n}. 
\end{align}
By the definition of $c_1$, we have 
\begin{align}
\frac{c_0 \sigma_Y^2 \norm{f}}{ q}
\absol{c_1\frac{\epsilon^X_n}{\epsilon^Y_n} - 1} 
\frac{1-e^{-M_Y(q) \epsilon^Y_n}}{\epsilon^Y_n}
\rightarrow \frac{c_0 \sigma_Y^2 \norm{f}}{ q}\times 0 \times M_Y(q)=0, 
~\text{ as } n\uparrow \infty. 
\end{align}
So we have \eqref{434}. 
By \eqref{432ac}, \eqref{433ab}, \eqref{434ac}, \eqref{434} and Lemma \ref{Lem402y}, we have 
\begin{align}
&\lim_{n\uparrow \infty} 
\frac{\sigma_X^2}{{(\epsilon^{X}_n)}^2}
\int_0^{\epsilon^X_n}R^{(q)}_{k_{T^+_0}Y^{(n)}}
f\rbra{\psi^{(n)}( v , v- \epsilon^X_n)}e^{-\Phi_{X^{(n)}}(q)v}dv  \\
=&c_0 n^Y_0 \sbra{\int_0^{T_0} e^{-qt}f(Y_t)dt; T^-_0=0} \\
=&n^U_0 \sbra{\int_0^{T_0} e^{-qt}f(U_t)dt; T^-_0=0}. \label{424}
\end{align}
\par
By \eqref{433z}, \eqref{429ab}, \eqref{431ab} and \eqref{424}, 
we obtain \eqref{423ab}. 
\par
{$ii$) We prove \eqref{424a} for $x\neq 0$. 
For $x<0$, we obtain \eqref{424a} by $i$) and the same argument as 
that of the proof of \cite[Theorem 8.4]{NobYan}. 
We now prove \eqref{424a} for $x>0$. 
We divide 
\begin{align}
R^{(q)}_{U}f(x)=
\bE^{U}_x\sbra{\int_0^{T^-_0}e^{-qt}f(U_t)dt}
+ \bE^{U}_x\sbra{\int_{T^-_0}^{T_0}e^{-qt}f(U_t)dt}
+\bE^{U}_x\sbra{\int_{T_0}^\infty e^{-qt}f(U_t)dt}
\end{align}
and we can divide $R^{(q)}_{U^{(n)}}f(x)$ similarly. 
By the definition of $U$ and \cite[Theorem 3.2]{Kyp}, we have the following: 
\begin{align}
\bE^{U}_x\sbra{\int_0^{T^-_0}e^{-qt}f(U_t)dt}
&=\bE^{X}_x\sbra{\int_0^{T^-_0}e^{-qt}f(X_t)dt}, \\
\bE^{U}_x\sbra{\int_{T^-_0}^{T_0}e^{-qt}f(U_t)dt}
&=\bE^{X}_x\sbra{e^{-qT^-_0}R^{(q)}_{k_{T^+_0}Y}f(\psi
(X_{T^-_0-}, X_{T^-_0})) }, \\
\bE^{U}_x\sbra{\int_{T_0}^\infty e^{-qt}f(U_t)dt}
&=\bE^{X}_x\sbra{e^{-qT^-_0}e^{\Phi_{Y} (q)\psi
(X_{T^-_0-}, X_{T^-_0})} }R^{(q)}_{U}f(0) ,
\end{align}
where we understand $\psi (0, 0)=0$. 
We have similar identities also for $U^{(n)}$. 
By the dominated convergence theorem, by the uniformly convergent coupling, 
by \cite[Lemma 8.3]{NobYan}, $i$) and by $\lim_{n\uparrow \infty}\Phi_{Y^{(n)}} =\Phi_Y$, 
it is sufficient to prove that 
\begin{align}
T^-_0(X^{(n)}) \rightarrow T^-_0 (X) \text{ and }
\psi^{(n)}(X^{(n)}_{T^-_0(X^{(n)})-}, X^{(n)}_{T^-_0(X^{(n)})}) \rightarrow 
\psi(X_{T^-_0(X)-}, X_{T^-_0(X)}) \label{463b}
\end{align}
hold as $n\uparrow \infty$ almost surely. 
\par
First, we prove \eqref{463b} on $A:=\{T^-_0(X)=\infty\} \cup \{T^-_0(X)<\infty, X_{T^-_0(X)}<0\}$. 
By the same argument as that of the proof of \cite[Theorem 8.4]{NobYan}, we have
	\begin{align}
	T^-_0 {(X)}= T^-_0 (X^{(n)}) \text{ for large } n \text{ on }A
	\end{align}   
and 
\begin{align}
\lim_{n \uparrow \infty} X^{(n)}_{T^-_0 (X^{(n)})-}=X_{T^-_0 (X)-} \text{ and }
\lim_{n \uparrow \infty} X^{(n)}_{T^-_0 (X^{(n)})}=X_{T^-_0 (X)} 
\text{ on }A.\label{463a}
\end{align}
By \eqref{463a} and the definition of $\psi^{(n)}$, we obtain 
\eqref{463b} on $A$. 
\par
Second, we prove \eqref{463b} on $A^c=\{T^-_0 (X)<\infty , X_{T^-_0 (X)}=0\}$. 
Let $\epsilon>0$ and let us argue on $A^c$. 
Set $I_{\epsilon}:=[T^-_0(X)-\epsilon , T^-_0(X)+\epsilon] $ and 
$\epsilon^\prime:=\rbra{\inf_{t\in[0, T^-_0(X)-\epsilon]}X_t} 
\land \absol{ \inf_{t\in I_\epsilon}X_t }$. 
Then there exists $N(\epsilon)>0$ such that 
for all $n>N(\epsilon)$, we have 
\begin{align}
\sup_{t\in[0, T^-_0 (X)+\epsilon]}\absol{X^{(n)}_t-X_t}< \epsilon^\prime. 
\label{463w}
\end{align}
By \eqref{463w}, \eqref{404} and the definition of $\psi^{(n)}$, for $n>N(\epsilon)$, we have
\begin{align}
&T^-_0(X)-\epsilon<T^-_0(X^{(n)})<T^-_0(X)+\epsilon, \label{464} \\
&\psi^{(n)}(X^{(n)}_{T^-_0(X^{(n)})-}, X^{(n)}_{T^-_0(X^{(n)})})< 2 
\rbra{l \lor c_1} 
\rbra{\sup_{t\in I_\epsilon}X_t-
\inf_{t\in I_\epsilon}X_t} . \label{465}
\end{align}
By \eqref{464} and \eqref{465}, we have \eqref{463b} on $A^c$. 
\par
The proof is therefor completed. 
}
}
\begin{Cor}\label{Cor403}
Under the same assumption of Theorem \ref{Thm404}, the process 
$(U^{(n)}, \bP^{U^{(n)}}_x)$ converges in distribution to 
$(U, \bP^{U}_x)$ for all $x \in \bR$. 
\end{Cor}
The proof of Corollary \ref{Cor403} can be obtained in the same way as that of 
\cite[Theorem 8.1 and 8.5]{NobYan} using scale functions of $U$ and $U^{(n)}$, 
so we omit it.

\section{Preliminary facts about duality}\label{duality}
In this section, we recall the definition of duality. 
\par
Let $\bT$, $Z$ and $m_Z$ be the same as those in Section \ref{Pre}. 
We assume that the process $(\hat{Z}, \bP^{\hat{Z}}_x)$ considered in this paper is a 
$\bT$-valued standard process with no negative jumps 
satisfying the following conditions: 
\begin{itemize}
\item[(C1)]{$(x, y) \mapsto \bE^{\hat{Z}}_x\sbra{ e^{-T_y}}>0$ is a ${\cal{B}}(\bT) \times \cal{B}(\bT)$-measurable function.} 
\label{iv}
\item[(C2)]{$\hat{Z}$ has a reference measure $m_Z$ on $\bT$.
}
\end{itemize}
\begin{Def}[See e.g., \cite{ChuWal}]
We say that $Z$ and $\hat{Z}$ are \textit{in duality} (relative to $m_Z$) if 
for $q >0$, non-negative measurable functions $f$ and $g$, we have
\begin{align}
\int_\bT f(x)R^{(q)}_Zg(x)m_Z(dx)=\int_\bT R^{(q)}_{\hat{Z}}f(x) g(x)m_Z(dx). 
\end{align}
\end{Def}
\begin{Thm}[See e.g., \cite{ChuWal} or \cite{Rev}]\label{Thm102a}
We suppose $Z$ and $\hat{Z}$ be  
in duality relative to $m_Z$. 
Then, for each $q > 0$, there exists a function $r_Z^{(q)}: \bT \times \bT \rightarrow 
[0, \infty)$ such that
\begin{enumerate}
\item{$r_Z^{(q)}$ is ${\cal{B}}(\bT)\times {\cal{B}}(\bT)$-measurable.}
\item{$x \mapsto r_Z^{(q)}(x, y) $ is $q$-excessive 
and finely continuous for each $y \in \bT$.}
\item{$y \mapsto r_Z^{(q)}(x, y) $ is $q$-coexcessive 
and cofinly continuous for each $x \in \bT$.}
\item{For all non-negative function $f$, 
\begin{align}
R_Z^{(q)}f(x)=\int_\bT f(y) r_Z^{(q)}(x, y) m_Z(dy), ~~~~~  
R_{\hat{Z}}^{(q)}f(y)=\int_\bT f(x) r_Z^{(q)}(x, y) m_Z(dx). \label{402x}
\end{align}
}
\end{enumerate}
\end{Thm}
By \cite[Proposition of Section $V$.1]{Rev}, 
if $Z$ and $\hat{Z}$ are in duality relative to $m_Z$, 
we normalize families of local times ${\{L^{Z, x}\}}_{x \in \bT} $ and 
${\{{{L}}^{\hat{Z}, x}\}}_{x \in \bT} $ 
which are 
{the same as} those in Section \ref{Pre} by 
\begin{align}
\bE^Z_x\sbra{\int_{0}^\infty e^{-qt} dL^{Z, y}_{t}}=r_Z^{(q)}(x , y), ~
\bE^{\hat{Z}}_y\sbra{\int_{0}^\infty e^{-qt} dL^{{\hat{Z}}, x}_{t}}=r_Z^{(q)}(x , y). \label{403x}
\end{align}
for all $q >0$. 
\begin{Lem}[{\cite[Lemma $4.16$]{GetSha}}]\label{Lem102a}
We assume that $Z$ and $\hat{Z}$ have the following conditions: 
\begin{itemize}
\item{$Z$ and $\hat{Z}$ are in duality relative to $m_Z$. }
\item{$Z$ and $\hat{Z}$ are recurrent processes. }
\item{$Z_{0}=Z_{T_x-}=x$, $n^Z_x$-a.s.. (This condition is equivalent to the counterpart of $n^{\hat{Z}}_x$.) }
\end{itemize}
Then 
we have 
\begin{align}
n^Z_x \sbra{~\cdot~} 
= n^{\hat{Z}}_x\sbra{\rho_x (~\cdot ~) }.  
\end{align}
\end{Lem} 
When $Z$ and $\hat{Z}$ are in duality, we always use 
the local times defined by \cite[Proposition of Section $V$.1]{Rev}. 
In other cases, we use the normalization of the local times in Section \ref{Pre}. 
We let scale functions ${\{ W_Z^{(q)}\}}_{q \geq 0}$ and 
${\{ W_{-\hat{Z}}^{(q)}\}}_{q \geq 0}$ be those in Section \ref{Pre}. 
Then we have the following lemma: 
\begin{Thm}[{\cite[Theorem 4.5]{Nob}}]\label{Lem206b}
If $Z$ and $\hat{Z}$ are in duality relative to $m_Z$, 
then we have
\begin{align}
W_Z^{(q)}(x, y) = W_{-{\hat{Z}}}^{(q)}(-y, -x),~~~~~x, y \in (b_0 , a_0 ) .   \label{123}
\end{align} 
If $\bT$ is open, then the converse is also true. 
\end{Thm}

\section{Duality problem of refracted processes}\label{SecDua}
In this section, we obtain the necessary and sufficient condition 
that the refracted processes $U$ and $\hat{U}$ are in duality 
in terms of an identity involving excursion measures 
and landing functions. 
\par
We assume that $\bT$ is an open set. 
Let $X$ and $Y$ be recurrent standard processes 
which are same as those in Section \ref{SecGen}. 
We assume that $0$ is irregular for itself for $X$ and $Y$ or 
$0$ is regular for itself for $X$ and $Y$.  
Let $\hat{X}$ and $\hat{Y}$ be $\bT_X$-valued and $\bT_Y$-valued standard 
processes with no negative jumps 
which satisfy the following conditions: 
\begin{itemize}
\item[(\^{B1})]{$(x, y) \rightarrow \bE^{\hat{X}}_x \sbra{e^{-T_y}}>0$ 
(resp. $(x, y) \rightarrow \bE^{\hat{Y}}_x \sbra{e^{-T_y}}>0$) is a 
${\cal{B}}(\bT_X) \times {\cal{B}}(\bT_X)$
(resp. ${\cal{B}}(\bT_Y) \times {\cal{B}}(\bT_Y)$)-measurable function. }
\item[(\^{B2})]
We assume that $\lim_{y\downarrow x}\bE^{\hat{X}}_y\sbra{e^{-T_x}}=1$ for all 
$x\in \bT_X\cap[0, \infty)$
(resp. $\lim_{y\downarrow x}\bE^{\hat{Y}}_y\sbra{e^{-T_x}}=1$ for all 
$x\in \bT_Y\cap (-\infty , 0)$). 
\item[(\^{B3})]
We assume that $\lim_{x \uparrow a_0}\bE^{\hat{X}}_x\sbra{ e^{- T^-_y}}=0$ for all $y\in\bT_X$
(resp. 
We assume that $\lim_{x \downarrow b_0}\bE^{\hat{Y}}_x\sbra{ e^{- T^+_y}}=0$ for all $y\in\bT_Y$). 
\item[(\^{B4})]{${\hat{X}}$ (resp. ${\hat{Y}}$) has a reference measure $m_X$ on $\bT_X$ 
(resp. $m_Y$ on $\bT_Y$).  }
\end{itemize}
In addition we assume the following conditions: 
\begin{itemize}
\item{
$X_0=X_{T_0-}=0$, $n^X_0$-a.s.. $Y_0=Y_{T_0-}=0$, $n^Y_0$-a.s.. }
\item{$X$ and $\hat{X}$(resp. $Y$ and $\hat{Y}$) are in duality 
relative to $m_X$ (resp. $m_Y$).  }
\end{itemize}
We take the local times ${\{L^{X, x}\}}_{x \in \bT_X}$, 
${\{L^{Y, x}\}}_{x \in \bT_Y}$, ${\{L^{\hat{X}, x}\}}_{x \in \bT_X}$ 
and ${\{L^{\hat{Y}, x}\}}_{x \in \bT_Y}$, the excursion measures 
${\{n^X_x\}}_{x \in \bT_X}$, ${\{n^Y_x\}}_{x \in \bT_Y}$, 
${\{n^{\hat{X}}_x\}}_{x \in \bT_X}$ and ${\{n^{\hat{Y}}_x\}}_{x \in \bT_Y}$, 
and the scale functions ${\{W_X^{(q)}\}}_{q \geq 0}$, 
${\{W_Y^{(q)}\}}_{q \geq 0}$, ${\{W_{-\hat{X}}^{(q)}\}}_{q \geq 0}$, 
and ${\{W_{-\hat{Y}}^{(q)}\}}_{q \geq 0}$ as those in Section \ref{duality}. 
As the landing functions, let $\psi:
(0, \infty) \times (-\infty , 0) \rightarrow (- \infty , 0)$ be 
a measurable function satisfying \eqref{301i} 
and $\phi: (- \infty , 0) \times (0, \infty) \rightarrow (0, \infty)$ be 
 a measurable function satisfying 
\begin{align}
n^Y_0 \sbra{1- e^{-T^+_0} \bE^X_{
\phi(Y_{ T^+_0-} , Y_{ T^+_0}) } 
\sbra{e^{-T_0}} ;0<T^+_0 < T_0} < \infty. \label{601}
\end{align}
Let $\bP^{U^0}_x$ and $n^U_0$ be those in Section \ref{SecGen}. 
By the excursion theory, we can construct 
a $\bT$-valued right continuous strong Markov processes 
$U$ from $n^U_0$ and ${\{\bP^{U^0}_x\}}_{x\in \bT \backslash \{0\}}$. 
Let $\hat{c}_0 \geq 0$ and {$\hat{c}_1 > 0$} be constants. 
We define the law of stopped process $\bP^{\hat{U}^0}_x$ 
for $x \neq 0$ and an excursion measure $n^{\hat{U}}_0$ away from $0$ 
by the following identities: 
\begin{align}
\bP^{{\hat{U}}^0}_x
\sbra{F({\hat{U}}^0) }
&= 
\begin{cases}
\bE^{{\hat{X}}^0}_x \sbra{ F({\hat{X}}^0)   },  ~~~~~~~~~~~~~~~~~~~~~~~~~~~~~~~~~~~~~~~~~~~~~~~  &x>0,   \\
\bE^{\hat{Y}}_x\sbra{
\bE^{{\hat{X}}^0}_{\phi ({\hat{Y}}_{T^+_0-}, {\hat{Y}_{T^+_0}})} \sbra{
F(\omega \circ {\hat{X}}^0)}\big{|}_{w={k_{T^+_0 } \hat{Y}}} ; T^+_0 \leq T_0 } 
 , ~&x< 0, 
\end{cases}
\\
n^{\hat{U}}_0 
\sbra{F(\hat{U}) }
&=\hat{c}_0 n^{\hat{X}}_0\sbra{ F( \hat{X}) ;T^+_0=0} \notag \\
&~+\hat{c}_1n^{\hat{Y}}_0\sbra{
\bE^{{\hat{X}}^0}_{\phi ({\hat{Y}}_{T^+_0-}, {\hat{Y}}_{T^+_0})} \sbra{
F(\omega \circ {\hat{X}}^0 )}\big{|}_{w=k_{T^+_0}\hat{Y}}
;0 < T^+_0 \leq T_0}    \label{306a}
\end{align}
for all positive measurable functional $F$. 
By the excursion theory, we can construct 
a $\bT$-valued right continuous strong Markov processes 
$\hat{U}$ from $n^{\hat{U}}_0$ together with ${\{\bP^{\hat{U}}_x\}}_{x \in \bT \backslash \{0\}}$.  
\par
We may and do assume $c_0=\hat{c}_0=\hat{c}_1=1$ without loss of generality.  
Let us explain the reason. 
We discuss positivity of $c_0$. 
By Lemma \ref{Lem102a}, the excursion measures $n^U_0$ and $n^{\hat{U}}_0$ 
need to satisfy $ n^U_0 \sbra{\cdot} 
= c_2n^{\hat{U}}_0\sbra{\rho_x (\cdot ) }$ for some constant $c_2>0$. 
This means that 
$n^U_0 \sbra{\cdot; T^-_0 = 0}= c_0 n^Y_0 \sbra{\cdot; T^-_0 = 0}
= c_2 \hat{c}_1n^{\hat{Y}}_0\sbra{\rho_x (\cdot ) ; T^+_0 = T_0}
= c_2n^{\hat{U}}_0\sbra{\rho_x (\cdot ) ; T^+_0 = T_0}$. 
So $c_0$ needs to be {equal to $c_2\hat{c}_1$} 
unless $n^{\hat{Y}}_0\sbra{\rho_x (\cdot ) ; T^+_0 = T_0}$ is the zero measure. 
When $n^{\hat{Y}}_0\sbra{\rho_x (\cdot ) ; T^+_0 = T_0}$ is the zero measure, 
so is $n^Y_0 \sbra{~\cdot~; T^-_0 = 0}$ by Lemma \ref{Lem102a}, 
which allows us to take $c_0 >0$. 
For the same reason, we may 
{assume that $1=c_2\hat{c}_0$. }
By changing the normalization of {$m_Y$, $n^Y_0$ and $n^{\hat{U}}_0$}, 
we may assume {$c_0=c_2=1$} without loss of generality, 
which yields {$\hat{c}_0 =\hat{c}_1=1$}. 
\par
We define $m_U = m_X|_{[0, \infty)} + m_Y|_{(-\infty, 0 )}$. Then we have the following 
theorem: 
\begin{Thm}\label{Thm501}
If $n^X_0$, $n^Y_0$, $\psi$ and $\phi$ 
satisfy 
\begin{align}
n^X_0 \sbra{ h ( X_{T^-_0-}, \psi (X_{T^-_0-} , X_{T^-_0})  ); 0<T^-_0<T_0}
=n^Y_0 \sbra{ h ( \phi (Y_{T^-_0} , Y_{T^-_0-}), Y_{T^-_0}  ); 0<T^-_0<T_0}\label{506aa}
\end{align}
for all non-negative measurable function $h$, or equivalently, 
\begin{align}
n^{\hat{X}}_0 \sbra{ h ( \hat{X}_{T^+_0}, \psi (\hat{X}_{T^+_0} , \hat{X}_{T^+_0-})  ); 0<T^+_0<T_0}
=n^{\hat{Y}}_0 \sbra{ h ( \phi (\hat{Y}_{T^+_0-} , \hat{Y}_{T^+_0}), \hat{Y}_{T^+_0-}  ); 0<T^+_0<T_0}
\label{605}
\end{align}
for all $h$, then $U$ and $\hat{U}$ are in duality relative to $m_U$. 
The converse is also true. 
\end{Thm}
\begin{Lem}\label{Lem302a}
If \eqref{506aa} is true, 
then we have
\begin{align}
n^U_0 \sbra{~ \cdot ~}\overset{d}{=} 
n^{\hat{U}}_0 \sbra{\rho_0(~\cdot~) }. \label{308b}
\end{align}
\end{Lem}
\Proof{
By \eqref{304a} and Lemma \ref{Lem102a}, 
for non-negative measurable functional $F$, we have 
\begin{align}
n^U_0 \sbra{F(U) }
&= n^X_0 \sbra{ F(X) ; T^-_0 =T_0}\notag \\
&+
n^X_0\sbra{\bE^{Y^0}_{\psi (X_{T^-_0-}, X_{T^-_0})} \sbra{
F(\omega \circ Y^0 )}\big{|}_{w=k_{T^-_0}X}
;0< T^-_0 < T_0}      \notag
\\
&+ n^Y_0\sbra{ F( Y ) ;T^-_0=0  }
\\
&= n^{\hat{X}}_0 \sbra{F(\rho_0 \hat{X}) ; T^+_0 = 0}
\notag\\
&+
n^{\hat{X}}_0\sbra{\bE^{Y^0}_{\psi({\hat{X}}_{T^+_0}, {\hat{X}}_{T^+_0 -})} \sbra{
F(\omega \circ Y^0 )}\big{|}_{w=k_{T_0}\rho_0\theta_{T^+_0}{\hat{X}}_t}
; 0 <T^+_0 < T_0 }      \notag
\\
&+n^{\hat{Y}}_0\sbra{ F( \rho_0{\hat{Y}} ) ;T^+_0= T_0}
 \label{311b}
\end{align}
By {\eqref{605}}, Lemma \ref{Lem102a} and Fubini's theorem, we have
\begin{align}
&n^{\hat{X}}_0\sbra{\bE^{Y^0}_{\psi({\hat{X}}_{T^+_0}, {\hat{X}}_{T^+_0 -})} \sbra{
F(\omega \circ Y^0 )}\big{|}_{w=k_{T_0}\rho_0\theta_{T^+_0}{\hat{X}}_t}
; 0 <T^+_0 < T_0 } \notag \\
&=n^{\hat{X}}_0\sbra{
\int \bP^{{\hat{X}}^0}_{{\hat{X}}_{T^+_0}}
\sbra{{\hat{X}}^0 \in d \omega}
\bE^{Y^0}_{\psi({\hat{X}}_{T^+_0}, {\hat{X}}_{T^+_0 -})} \sbra{
F( k_{T_0}\rho_0\omega \circ Y^0)}
;0<T^+_0 < T_0 } \\
&=n^{\hat{Y}}_0\sbra{
\int \bP^{{\hat{X}}^0}_{\phi (\hat{Y}_{T^+_0-}, \hat{Y}_{T^+_0} )}
\sbra{{\hat{X}}^0 \in d \omega}
\bE^{Y^0}_{\hat{Y}_{T^+_0-}} \sbra{
F( k_{T_0}\rho_0\omega \circ Y^0)}
;0<T^+_0 < T_0 } \\
&=
n^{\hat{Y}}_0\sbra{
\bE^{Y^0}_{{\hat{Y}}_{T^+_0 -}} \sbra{\int 
F(\omega \circ Y^0 )
\bP^{{\hat{X}}^0}_{y}
\sbra{k_{T_0}\rho_0{\hat{X}} \in d \omega}}
\Big{|}_{y = \phi({\hat{Y}}_{T^+_0-}, {\hat{Y}}_{T^+_0 })}
;0<T^+_0 < T_0}  . 
\label{312b}
\end{align}
By the strong Markov property and Lemma \ref{Lem102a}, we have 
\begin{align}
\eqref{312b}
&=
n^{Y}_0\sbra{
\bE^{Y^0}_{Y_{T^-_0 }} \sbra{\int 
F(\omega \circ Y^0 )
\bP^{{\hat{X}}^0}_{y}
\sbra{k_{T_0}\rho_0{\hat{X}} \in d \omega}}
\Big{|}_{y = \phi(Y_{T^-_0}, Y_{T^-_0- })}
;0<T^-_0 < T_0}   \\
&=
n^{Y}_0\sbra{\int 
F(\omega \circ \theta_{T^-_0}Y )
\bP^{{\hat{X}}^0}_{ \phi ({Y}_{T^-_0}, {Y}_{T^-_0 -})}
\sbra{k_{T_0}\rho_0{\hat{X}} \in d \omega}
;0<T^-_0 < T_0}  
\\
&=
n^{Y}_0\sbra{
\bE^{{\hat{X}}^0}_{\phi({Y}_{T^-_0}, {Y}_{T^-_0 -})}
\sbra{
F(k_{T_0}\rho_0{\hat{X}} \circ \omega^\prime)
}\Big{|}_{\omega^\prime = \theta_{T^-_0}Y}
;0<T^-_0 < T_0}
\\  
&=
n^{\hat{Y}}_0\sbra{
\bE^{{\hat{X}}^0}_{\phi ({\hat{Y}}_{T^+_0-}, {\hat{Y}}_{T^+_0 })}
\sbra{
F(\rho_0( \omega  \circ {{\hat{X}}}^0))
}\Big{|}_{\omega = k_{T^+_0}\hat{Y}}
;0<T^+_0 <T_0 }. \label{316b}
\end{align}
By \eqref{316b}, we have 
\begin{align}
\eqref{311b}
&=n^{\hat{Y}}_0\sbra{ F( \rho_0{\hat{Y}} ) ;T^+_0= T_0} \notag \\
&+n^{\hat{Y}}_0\sbra{
\bE^{{\hat{X}}^0}_{\phi ({\hat{Y}}_{T^+_0-}, {\hat{Y}}_{T^+_0 })}
\sbra{
F(\rho_0( \omega  \circ {\hat{X}}^0))
}\Big{|}_{\omega = k_{T^+_0}\hat{Y}}
;0<T^+_0 <T_0 } \notag \\
&+n^{\hat{X}}_0\sbra{ F( \rho_0{\hat{X}}) ;T^+_0=0}  \\
&=n^{\hat{U}}_0 \sbra{F ( \rho_0 \hat{U}) }
\end{align}
So we obtain \eqref{308b}. 
}
\begin{Lem}\label{Lem503}
For all $q > 0$ and $x\in\bT$, 
{the measure $R^{(q)}_U1_{(\cdot)}(x)$ is absolutely continuous 
with respect to $m_U ( \cdot)$}. 
\end{Lem}
\Proof{
Let $A$ be a set in ${\cal{B}(\bT)}$ which satisfies 
$m_X (A\cap [0 ,\infty))= 0$ and $m_Y (A\cap (-\infty, 0))= 0$. 
It is sufficient to prove that $\bE^U_0 \sbra{\int_0^\infty e^{-qt}1_A(U_t)dt}=0$. 
By the compensation theorem of excursion point processes, 
we have 
\begin{align}
&qn^U_0 \sbra{1-e^{-qT_0}}\bE^U_0 \sbra{\int_0^\infty e^{-qt}1_A (U_t) dt} \\
&={n^U_0\sbra{\int_0^{T_0} e^{-qt}1_A(U_t)dt }}\\
&={n^X_0\sbra{\int_0^{T^-_0} e^{-qt}1_A(X_t)dt}  
+n^X_0\sbra{\bE^Y_{J_X}
\sbra{\int_0^{T^+_0} e^{-qt}1_A(Y_t)dt}}
+n^Y_0\sbra{\int_0^{T_0 } e^{-qt} 1_A(Y_t)dt ; T^-_0=0}}. 
\end{align}
By the assumption of $A$, we have 
\begin{align}
n^X_0 \sbra{\int_0^{T^-_0} e^{-qt}1_A(X_t)dt}
=qn^X_0\sbra{1-e^{-qT_0}}\bE^X_0 
\sbra{\int_0^\infty e^{-qt}1_{A\cap [0, \infty)}(X_t) dt}=0, 
\end{align} 
\begin{align}
\bE^Y_{J_X}
\sbra{\int_0^{T^+_0} e^{-qt}1_A(Y_t)dt}\leq 
\bE^Y_{J_X}
\sbra{\int_0^{\infty} e^{-qt}1_{A\cap(-\infty , 0)}(Y_t)dt}=0 .
\end{align}
and 
\begin{align}
n^Y_0\sbra{\int_0^{T_0 } e^{-qt} 1_A(Y_t)dt ; T^-_0=0}
\leq qn^Y_0\sbra{1-e^{-qT_0}}\bE^Y_0 
\sbra{\int_0^\infty e^{-qt}1_{A\cap (-\infty, 0)}(Y_t) dt}=0
\end{align}
So we obtain $\bE^U_0 \sbra{ \int_0^\infty e^{-qt}1_A(U_t)dt}=0$. 
}
We want to find suitable normalization of local times of $U$. 
By \cite[Theorem $18.4$]{GemHor}, 
we let local times 
${\{L^{U,x \prime}\}}_{x \in \bT\backslash \{0\}}$ of $U$ be those in Section \ref{Pre}. 
We set $n^{U\prime}_0 =n^U_0$ and let $n^{U\prime}_x$ 
for $x\in \bT\backslash \{0\}$ be the excursion measure associated to 
$L^{U, x \prime}$. 
Then there exists the positive function $c(x)$ such that $c(0)=1$ 
(by the definition of $U^0$) and 
for all non-negative functional $F$: 
\begin{align}
n^{U \prime}_x \sbra{F({\{U_t \}}_{t < T^-_0})} 
&=c(x)n^X_x \sbra{F({\{X_t \}}_{t < T^-_0})},~~~x\in\bT\cap[0, \infty), \label{325a}\\
n^{U\prime}_x \sbra{F({\{U_t \}}_{t < T^+_0})} 
&=c(x)n^Y_x \sbra{F({\{Y_t \}}_{t < T^+_0})},~~~x\in\bT\cap(-\infty , 0].\label{326a}
\end{align} 
Then we have $c(x)=1$ $m_U$-a.e. 
Indeed, 
for all $q>0$, $x, y \in \bT\cap [0, \infty)$ 
and non-negative measurable function $f$, we have 
\begin{align}
\bE^U_x \sbra{ \int_{0}^{T^-_0}e^{-qt} dL^{U, y \prime}_t }  
=\frac{1}{c(y)}
 \bE^X_x \sbra{ \int_{0}^{T^-_0}e^{-qt} dL^{X, y}_t } 
\end{align}
and
\begin{align}
\int_{\bT\cap [0, \infty)} f(y) \bE^U_x \sbra{ \int_{0}^{T^-_0}e^{-qt} dL^{U, y \prime}_t }  m_U (dy)
&=\bE^U_x \sbra{\int_0^{T^-_0} e^{-qt} f(U_t) dt } \\
&=\bE^X_x \sbra{\int_0^{T^-_0} e^{-qt} f(X_t) dt } \\
&=\int_{\bT\cap [0, \infty)} f(y) \bE^X_x \sbra{ \int_{0}^{T^-_0}e^{-qt} dL^{X, y}_t }  m_U (dy).
\end{align}
So $c(x)=1$ on $\bT\cap [0, \infty)$ $m_U$-a.e. 
Similarly, $c(x)=1$ on $\bT\cap (- \infty, 0)$ $m_U$-a.e. 
We now set $L^{U, x}=c(x)L^{U, x, \prime}$ 
and $n^U_x = \frac{1}{c(x)}n^{U, \prime}_x$. 
This local times satisfy \eqref{201f} and \eqref{202f} since $c(x)=1$ $m_U$-a.e. 
\par
In the same way, 
let the excursion measures ${\{n^{\hat{U}}_x\}}_{x\in\bT}$ of $\hat{U}$ be those in 
Section \ref{Pre} satisfying the following conditions; 
\begin{align}
n^{\hat{U}}_x \sbra{F({\{{\hat{U}}_t \}}_{t < T^-_0})} 
&=n^{\hat{X}}_x \sbra{F({\{{\hat{X}}_t \}}_{t < T^-_0})},~~~x\in\bT \cap [0, \infty), \\
n^{\hat{U}}_x \sbra{F({\{{\hat{U}}_t \}}_{t < T^+_0})} 
&=n^{\hat{Y}}_x \sbra{F({\{{\hat{Y}}_t \}}_{t < T^+_0})},~~~x\in\bT\cap(-\infty, 0]. 
\end{align}
We let the scale functions $\cbra{W_U^{(q)}}_{q \geq 0}$ and 
$\cbra{W_{-\hat{U}}^{(q)}}_{q \geq 0}$ be those in \eqref{216}. 
\Proof[Proof of Theorem \ref{Thm501}]{
Let us assume that we have \eqref{506aa} for all non-negative measurable function $h$. 
By Theorem \ref{Lem206b} and Lemma \ref{Lem503}, it is sufficient to prove that 
\begin{align}
W_U^{(q)}(x, y) =W_{-\hat{U}}^{(q)}(-y, -x),\label{331a}
\end{align}
for $q \geq 0$ and $x, y \in \bT$. 
For $0\leq y < x$, we have 
\begin{align}
&~~~~~~W_U^{(q)}(x, y)
={n^U_y\sbra{e^{-q T^+_x} ;T^+_x< \infty}}^{-1}
={n^X_y\sbra{e^{-q T^+_x} ;T^+_x< \infty}}^{-1}
=W_X^{(q)}(x, y)\\
&=W_{-\hat{X}}^{(q)}(-y, -x)
={n^{-\hat{X}}_{-x}\sbra{e^{-qT^+_{-y}};T^+_{-y} < \infty }}^{-1}
={n^{-\hat{U}}_{-x}\sbra{e^{-qT^+_{-y}};T^+_{-y} < \infty }}^{-1}
=W_{-\hat{U}}^{(q)}(-y, -x) \label{333a}
\end{align}
by the definitions of $n^U_y$, $n^{-\hat{U}}_{-x}$ and Theorem \ref{Lem206b}. 
Similarly, for $y< x \leq 0$, we have
\begin{align}
W_U^{(q)}(x, y)=W_Y^{(q)}(x, y)=W_{-\hat{Y}}^{(q)}(-y, -x)
=W_{-\hat{U}}^{(q)}(-y, -x). \label{334a}
\end{align}
When $y< 0< x$, by \eqref{122aa}, {\eqref{216} and \eqref{103},} we have
\begin{align}
W_U^{(q)}(x, y) = 
	W_U^{(q)}(0, y)W_U^{(q)}(x, 0)
n^{U}_0\sbra{ 1 - e^{-qT_0}1_{\{T^-_y = \infty, T^+_x = \infty  \}}}
 \label{335a}
\end{align}
and 
\begin{align}
W_{-\hat{U}}^{(q)}(-y, -x)=
W_{-\hat{U}}^{(q)}(-y, 0)W_{-\hat{U}}^{(q)}(0, -x)
n^{-\hat{U}}_0\sbra{ 1 - e^{-qT_0}1_{\{T^-_y = \infty, T^+_x = \infty  \}}}. 
\label{535aa}
\end{align}
By Lemma \ref{Lem302a}, \eqref{333a}, \eqref{334a}, 
\eqref{335a} and \eqref{535aa}, we obtain \eqref{331a}. \par
We assume that $U$ and $\hat{U}$ are in duality relative to $m_U$. 
By Lemma \ref{Lem102a} and the definitions of $n^U_0$ and $n^{\hat{U}}_0$, 
we have \eqref{308b}. 
We have 
\begin{align}
n^X_0 \sbra{ h ( X_{T^-_0-}, \psi (X_{T^-_0-} , X_{T^-_0})  ); 0<T^-_0<T_0}
=n^U_0 \sbra{ h ( U_{T^-_0-}, U_{T^-_0}  ); 0<T^-_0<T_0}\label{537aa}
\end{align}
and
\begin{align}
n^Y_0 \sbra{ h ( \phi (Y_{T^-_0} , Y_{T^-_0-}), Y_{T^-_0}  ); 0<T^-_0<T_0}
&=n^{\hat{Y}}_0 \sbra{ h ( \phi ({\hat{Y}}_{T^+_0-} , Y_{T^+_0}), {\hat{Y}}_{T^+_0-}  ); 0<T^+_0<T_0}\\
&=n^{\hat{U}}\sbra{ h ( {\hat{U}}_{T^+_0}, {\hat{U}}_{T^+_0-}  ); 0<T^+_0<T_0}.\label{539aa}
\end{align}
By \eqref{308b}, \eqref{537aa} and \eqref{539aa}, we obtain \eqref{506aa}. 
The proof is completed. 
}

\section{An example of the duality problem}\label{Exa}
In this section, we construct refracted processes in duality 
from spectrally negative stable processes. 
\par
Let $X$ be a spectrally negative strictly $\alpha$-stable process 
whose L\'evy measure is 
\begin{align}
\Pi_X(dx) =c_X 1_{\{x < 0\}}{\absol{x}}^{-\alpha -1}dx
\end{align}
for a constant $c_X > 0 $, 
and 
$Y$ be a spectrally negative strictly $\beta$-stable process 
whose L\'evy measure is 
\begin{align}
\Pi_Y(dx) =c_Y 1_{\{x < 0\}}{\absol{x}}^{-\beta -1}dx 
\end{align}
where $c_Y > 0$. 
Then it is known that
\begin{align}
\hat{X}~ (\text{under } ~\bP^{\hat{X}}_x) \overset{d}{=} -X ~(\text{under } ~\bP^{X}_{-x})
\end{align}
and 
\begin{align}
\hat{Y}~ (\text{under } ~\bP^{\hat{Y}}_x) \overset{d}{=} -Y ~(\text{under } ~\bP^{Y}_{-x}).
\end{align}
We set reference measure $m_X(dx)$ as $\frac{\alpha-1}{c_X}dx$ 
and 
reference measure $m_Y(dx)$ as $\frac{\beta-1}{c_Y}dx$. 
Let $n^X_0$ and $n^Y_0$ be those in Section \ref{SecDua}. 
We want to find suitable landing functions such that $U$ 
and $\hat{U}$ are in duality. 
So we need to find $\psi$ and $\phi$ satisfying \eqref{506aa}. 
\begin{Prop}
Suppose $\alpha> \beta$. 
We let $\psi(x, y)=y{(x-y)}^{\frac{\alpha-1}{\beta-1}-1}$ and 
$\phi(x, y) =y {(y-x)}^{\frac{\beta-1}{\alpha-1} - 1}$. 
Then $U$ constructed from $X$, $Y$, $\psi$ and $c_0=0$ and 
$\hat{U}$ constructed from $\hat{X}$, $\hat{Y}$, $\hat{\psi}$ and $c_0=0$ 
are well-defined and in duality relative to $m_U$. 
\end{Prop}
\Proof{
Let us prove \eqref{506aa}. 
By \cite[Theorem 3.3]{NobYan}, we have 
\begin{align}
&n^{X}_0\sbra{
		h(X_{T^-_0 - }, \psi(X_{T^{-}_{0}-}, X_{T^{-}_{0}}) ); 0< T^{-}_{0} < T_0}
	\notag \\ 
		&=\frac{\alpha-1}{c_X}\int_0^\infty dv \int_{(- \infty , 0)}h(v, \psi(v, u)) {\Pi}_{X}(du -v) \\
&=(\alpha-1) \int_0^\infty dv \int_0^\infty h(v, \psi(v, -u)) 
	{(u+v)}^{-1-\alpha}du \label{505j}
\end{align}
and
\begin{align}
&n^{Y}_0\sbra{
		h( \phi(Y_{T^-_0}, Y_{T^-_0 -}), Y_{T^-_0} ); 0< T^{-}_{0} < T_0}
	\notag \\ 
		&=\frac{\beta-1}{c_Y}\int_0^\infty dv \int_{(- \infty, 0)}h(\phi(u, v), u) {\Pi}_{Y}(du -v) \\
&=(\beta-1) \int_0^\infty dv \int_0^\infty h(\phi(-u, v), -u) 
	{(u+v)}^{-1-\beta}du .\label{507j}
\end{align} 
We set $s=\frac{u}{u+v}$, $t = u+v$, $t_1 = t^{-\alpha+1}$ and $t_2=t^{-\beta+1}$. 
Then we have $u=st$, $v=t(1-s)$ and 
$\absol{\frac{\partial u}{\partial s}\frac{\partial v}{\partial t} 
-\frac{\partial u}{\partial t}\frac{\partial v}{\partial s}}=t$. 
So we have 
\begin{align}
\eqref{505j}
&=(\alpha-1)
\int_0^1 ds \int_0^\infty 
h(t(1-s) , \psi(t(1-s) , -st)) t^{-\alpha} dt \\
&=\int_0^1 ds \int_0^\infty 
h({t_1}^{-\frac{1}{\alpha-1}}(1-s) , \psi({t_1}^{-\frac{1}{\alpha-1}}(1-s) , -s{t_1}^{-\frac{1}{\alpha-1}}))  dt_1 \label{511j}
\end{align}
and 
\begin{align}
\eqref{507j}
&=(\beta-1)\int_0^1 ds \int_0^\infty 
h( \phi(-st, t(1-s) ) , -st ) t^{-\beta} dt \\
&=\int_0^1 ds \int_0^\infty 
h( \phi(-s{t_2}^{-\frac{1}{\beta-1}}, {t_2}^{-\frac{1}{\beta-1}}(1-s) ) , -s{t_2}^{-\frac{1}{\beta-1}} )  dt_2.  \label{513j}
\end{align}
Since $\psi(x, y)=y{(x-y)}^{\frac{\alpha-1}{\beta-1}-1}$ and 
$\phi(x, y) =y {(y-x)}^{\frac{\beta-1}{\alpha-1} - 1}$, we have 
\begin{align}
\psi({t}^{-\frac{1}{\alpha-1}}(1-s) , -s{t}^{-\frac{1}{\alpha-1}})=-s{t}^{-\frac{1}{\beta-1}},
~~~~~~s, t>0\label{514j}
\end{align}
and 
\begin{align}
\phi(-s{t}^{-\frac{1}{\beta-1}}, {t}^{-\frac{1}{\beta-1}}(1-s) )
={t}^{-\frac{1}{\alpha-1}}(1-s),~~~~~~s, t>0 . \label{515j}
\end{align}
By \eqref{511j}, \eqref{513j}, \eqref{514j} and \eqref{515j}, we obtain \eqref{506aa}. 
\par
Let us prove \eqref{301i} and \eqref{601}. 
Let $\Phi_X$ and $\Phi_Y$ be those in Section \ref{App}. 
By \cite[Theorem 3.2]{Kyp}, we have 
\begin{align}
&n^X_0 \sbra{1- e^{-T^-_0} \bE^Y_{
\psi(X_{ T^-_0-} , X_{ T^-_0}) } 
\sbra{e^{-T_0}} ;0<T^-_0 < T_0}  \\
=&n^X_0 \sbra{1- e^{-T^-_0} 
e^{\Phi_Y(1){\psi(X_{ T^-_0-} , X_{ T^-_0}) }} 
 ;0<T^-_0 < T_0} \\
\leq&
n^X_0 \sbra{1- e^{-T^-_0} 
e^{\Phi_Y(1){ X_{ T^-_0} }} 
 ;0<T^-_0 < T_0, X_{T^-_0-}- X_{T^-_0} \leq 1} \label{717}\\
&+n^X_0 \sbra{1- e^{-T^-_0} 
e^{\Phi_Y(1){\psi(X_{ T^-_0-} , X_{ T^-_0}) }} 
 ;0<T^-_0 < T_0, X_{T^-_0-}- X_{T^-_0} > 1} \label{718}
\end{align}
where the inequality \eqref{717} uses $\alpha > \beta$. 
There is a constant $q\geq 1$ such that $\Phi_X(q)\geq\Phi_Y (1)$. 
By \cite[Theorem 3.2]{Kyp} and the strong Markov property, we have 
\begin{align}
\eqref{717}
\leq n^X_0\sbra{  1-e^{-qT^-_0}e^{\Phi_X(q) X_{T^-_0}}}
=n^X_0\sbra{1-e^{-qT_0}}<\infty. \label{719}
\end{align}
By the property of excursion measures, we have 
\begin{align}
\eqref{718}\leq n^X_0\sbra{X_{T^-_0-}-X_{T^-_0}>1}<\infty. \label{720}
\end{align}
By \eqref{719} and \eqref{720}, we obtain \eqref{301i}. 
Since we have \eqref{506aa} and \eqref{301i}, 
in the same way as the proof of Lemma \ref{Lem302a}, we obtain \eqref{601}. 
So the refracted processes $U$ and $\hat{U}$ are well-defined. 
By \eqref{506aa} and Theorem \ref{Thm501}, 
the proof is completed. 
}

\section*{Acknowledgments}
I would like to express my deepest gratitude to 
Professor V\'ictor Rivero and my supervisor Professor Kouji Yano 
for comments and improvements. 
Especially, Professor Kouji Yano gave me a lot of advice. 
The author was supported by JSPS-MAEDI Sakura program. 


\end{document}